\documentclass[10pt]{article}
\usepackage[left=0.8in,top=0.8in,right=0.8in,bottom=0.8in]{geometry}
\usepackage{graphicx}
\usepackage{amsmath}
\usepackage{amssymb}
\usepackage{amsthm}
\usepackage{bm}
\usepackage{stmaryrd}
\usepackage{mathrsfs}
\usepackage{array}
\usepackage{multirow}
\usepackage{url}

\SetSymbolFont{stmry}{bold}{U}{stmry}{m}{n}

\newcolumntype{P}[1]{>{\centering\arraybackslash}p{#1}}

\begin{document}
\title{A note on the accuracy of the generalized-$\alpha$ scheme for the incompressible Navier-Stokes equations}
\author{Ju Liu, Ingrid S. Lan, Oguz Z. Tikenogullari, and Alison L. Marsden\\
\textit{\small Department of Pediatrics (Cardiology), Department of Bioengineering,}\\
\textit{\small Department of Mechanical Engineering,}\\
\textit{\small and Institute for Computational and Mathematical Engineering,}\\
\textit{\small Stanford University, Clark Center E1.3, 318 Campus Drive, Stanford, CA 94305, USA}\\
\small \textit{E-mail address:} \{liuju,ingridl,oguzziya,amarsden\}@stanford.edu}

\date{}
\maketitle

\section*{Abstract}
We investigate the temporal accuracy of two generalized-$\alpha$ schemes for the incompressible Navier-Stokes equations. In a popular approach, the pressure is collocated at the time step $t_{n+1}$ while the remainder of the Navier-Stokes equations is discretized following the rule of the generalized-$\alpha$ scheme. This scheme has been claimed to be second-order accurate in time. We developed a suite of numerical codes using inf-sup stable higher-order non-uniform rational B-spline (NURBS) elements for spatial discretization. In doing so, we are able to achieve high spatial accuracy and to perform temporal refinement without introducing stabilization terms, which degenerate at small time steps. Numerical evidence suggests that only first-order accuracy is achieved, at least for the pressure, in this aforesaid temporal discretization approach. Evaluating the pressure at the intermediate time step recovers second-order accuracy, and the numerical implementation becomes simpler. The second approach is recommended as the generalized-$\alpha$ scheme of choice when integrating the incompressible Navier-Stokes equations.


\vspace{5mm}

\noindent \textbf{Keywords:} Generalized-$\alpha$ scheme, Incompressible Navier-Stokes equations, Temporal accuracy

\section{Introduction}
\label{sec:introduction}
The generalized-$\alpha$ method was initially proposed for structural dynamics \cite{Chung1993} and possesses the desirable numerical properties of implicit schemes, as noted by Hilber and Hughes \cite{Hilber1978}. The method was later extended to address first-order systems, with an important example being the Navier-Stokes equations \cite{Jansen2000}. Interestingly, when used in conjunction with the first-order structural dynamics, the generalized-$\alpha$ method has been found not to suffer from the `overshoot' phenomenon \cite{Kadapa2017}, making it an ideal time integration algorithm for low-frequency, inertial type problems. Over the years, this temporal scheme has gained popularity as an implicit scheme for fluid and fluid-structure interaction (FSI) problems. When it was originally introduced to computational fluid dynamics (CFD) problems, Jansen et al. considered the compressible Navier-Stokes equations written in the pressure primitive variables, which were evaluated at the intermediate time step \cite{Jansen2000}. Later, when dealing with the incompressible Navier-Stokes equations, the time-stepping scheme was used in a different manner: the velocity was evaluated at the intermediate time step following the rule of the generalized-$\alpha$ scheme, while the pressure was collocated at the time step $t_{n+1}$. This choice of temporal discretization has in fact become quite popular in both CFD and FSI simulations (see e.g. \cite{svSolver-github-repo,Bazilevs2007a,Bazilevs2006,Bazilevs2013b,Casquero2016,Casquero2020, Figueroa2006,Gamnitzer2010,Kang2012,Joshi2019,Peters2019, Seo2019,Vignon-Clementel2006,Yan2018,Zhu2020}). In geometrically multiscale coupled problems, this dichotomy of time discretization treatment further leads to a competing question: whether one should collocate the boundary traction at the time step $t_{n+1}$ \cite{Moghadam2013} or at the intermediate time step \cite{Kim2009}. We mention that this approach is also adopted when applying the $\theta$-scheme for the incompressible Navier-Stokes equations \cite{Codina2007,Jodlbauer2019}, although the $\theta$-scheme is beyond the scope of this study. Indeed, the pressure is often viewed as a Lagrange multiplier and should therefore ``not be subjected to a time integration scheme" \cite{Saksono2007} and ``requires an implicit treatment" \cite[p.284]{Rannacher2008}. While this may explain the widespread adoption of the above mentioned scheme, we note that there exist a few works evaluating the pressure at the intermediate time step \cite{Gravemeier2011,Kadapa2020,Liu2018,Rossi2016}. In our view, treating the pressure differently in time is unnatural and inefficient for the following reasons. First, given that the incompressible Navier-Stokes equations constitute a limiting case of their compressible counterparts, this approach is evidently inconsistent with the scheme proposed in \cite{Jansen2000}. Second, in the implementation of CFD algorithms, variables are typically stored in an `interlaced' pattern, which is cache-friendly for high-performance computing \cite{Gropp2000a}. Using two distinct solution updating strategies in the time integration algorithm necessitates separate extraction of the velocity and pressure degrees-of-freedom from the interlaced solution vector, thereby reducing the overall algorithmic efficiency. In addition to the inconsistency and inefficiency, we may also reasonably expect a loss of numerical accuracy in the widely adopted approach.

In this short communication, we compare the temporal accuracy of two generalized-$\alpha$ schemes. The first, which we denote as Scheme-1, is the popular choice in which the pressure is collocated at the time step $t_{n+1}$; the second, denoted as Scheme-2, evaluates the pressure at the intermediate time step. We present the accuracy of quantities of interest in two benchmark cases. The Navier-Stokes equations are discretized in space by inf-sup stable NURBS element pairs. Our rationale behind this choice is twofold. We want to avoid invoking the stabilized formulation, which can degenerate for small time steps when used without special treatment \cite{Hsu2010} and is therefore not robust for temporal convergence studies. Furthermore, the robustness and exact geometric representation made possible by the NURBS-based isogeometric analysis technique allow us to eliminate the possibility of confounding errors from the spatial discretization and thereby focus on the temporal accuracy.

\section{Governing Equations and the Spatiotemporal Discretization}
\label{sec:continuum_problem_and_discretization}
\subsection{Strong-Form Problem}
Let $\Omega \subset \mathbb R^d$ be a fixed bounded open set with sufficiently smooth boundary $\Gamma := \partial \Omega$, where $d$ represents the number of spatial dimensions. The time interval is denoted $(0,T) \subset \mathbb R$ with $T>0$. The governing equations for the incompressible flow of a Newtonian fluid can be stated as
\begin{align}
\label{eq:ns_mom}
\bm 0 &= \rho \frac{\partial \bm v}{\partial t}  + \rho  \bm v \cdot \nabla \bm v - \nabla \cdot \bm \sigma - \rho \bm f, && \mbox{ in } \Omega \times (0,T), \\
\label{eq:ns_mass}
0 &= \nabla \cdot \bm v, && \mbox{ in } \Omega \times (0, T),
\end{align}
wherein
\begin{align*}
\bm \sigma := 2\mu \bm \varepsilon(\bm v) - p \bm I, \quad \bm \varepsilon(\bm v) := \frac12 \left( \nabla \bm v + \nabla \bm v^T \right).
\end{align*}
In the above, $\rho$ is the fluid density, $\bm v$ is the velocity field, $\bm \sigma$ is the Cauchy stress, $\bm f$ is the body force, $\mu$ is the dynamic viscosity, $\bm \varepsilon$ is the rate-of-strain tensor, $p$ is the pressure, and $\bm I$ is the second-order identity tensor. Both $\rho$ and $\mu$ are assumed constant, and the kinematic viscosity is defined as $\nu := \mu / \rho$. Given a divergence-free velocity field $\bm v_0$, the initial condition is given by
\begin{align}
\label{eq:initial_condition}
\bm v(\cdot, 0) = \bm v_0(\cdot), && \mbox{ in } \bar{\Omega}.
\end{align}
The boundary $\Gamma$ can be partitioned into two non-overlapping subdivisions, that is, $\overline{\Gamma} = \overline{\Gamma_g \cup \Gamma_h}$ and $\emptyset = \Gamma_g \cap \Gamma_h$. The subscripts $g$ and $h$ indicate the Dirichlet and Neumann partitions. The unit outward normal vector to $\Gamma$ is denoted as $\bm n$. Given the Dirichlet values $\bm g$ and the boundary traction $\bm h$, the initial-boundary value problem is to solve for $\bm v$ and $p$ that satisfy \eqref{eq:ns_mom}, \eqref{eq:ns_mass}, \eqref{eq:initial_condition}, and
\begin{align}
\label{eq:dirichlet_bc}
\bm v = \bm g  \mbox{ on } \Gamma_{g} \times (0,T), \quad \bm \sigma \cdot \bm n = \bm h  \mbox{ on } \Gamma_{h} \times (0,T).
\end{align}

\subsection{Semi-Discrete Formulation}
\label{subsec:semi_discrete_formulation}
\subsubsection{Spline Spaces on the Parametric Domain}
We start by reviewing the construction of B-splines and NURBS basis functions. Given the polynomial degree $\mathsf p$ and the dimensionality of the B-spline space $\mathsf n$, the knot vector can be represented by $\Xi := \left\{\xi_1, \xi_2 \cdots, \xi_{\mathsf n + \mathsf p + 1} \right\}$, wherein $0=\xi_1 \leq \xi_2 \leq \cdots \leq \xi_{\mathsf n + \mathsf p + 1}=1$.
The B-spline basis functions of degree $\mathsf p$, denoted as $\mathsf N_{i}^{\mathsf p}(\xi)$, for $i=1,\cdots, \mathsf n$, can then be defined recursively from the knot vector using the Cox-de Boor recursion formula \cite{Hughes2005}. Given a set of weights $\{\mathsf w_1, \mathsf w_2, \cdots , \mathsf w_{\mathsf n} \}$, the NURBS basis functions of degree $\mathsf p$ can be defined as
\begin{align*}
\mathsf R_i^{\mathsf p}(\xi) := \mathsf w_i \mathsf N^{\mathsf p}_i(\xi)/\mathsf W(\xi), \quad \mathsf W(\xi) := \sum\limits_{j=1}^{\mathsf n}\mathsf w_j \mathsf N^{\mathsf p}_j(\xi).
\end{align*}
Importantly, the knots can also be represented with two vectors, one of the unique knots $\{\zeta_1, \zeta_2, \cdots, \zeta_{\mathsf m} \}$ and another of the corresponding knot multiplicities $\{r_1, r_2, \cdots, r_{\mathsf m}\}$. As is standard in the literature of computer-aided design, we consider open knot vectors in this work, meaning that $r_1 = r_{\mathsf m} = \mathsf p+1$. We further assume that $r_i \leq \mathsf p$ for $i=2,\cdots, \mathsf m-1$. Across any given knot $\zeta_i$, the B-spline basis functions have $\alpha_i := \mathsf p - r_i$ continuous derivatives. The vector $\bm \alpha := \{ \alpha_1, \alpha_2, \cdots, \alpha_{\mathsf m-1}, \alpha_m \} = \{ -1, \alpha_2, \cdots, \alpha_{\mathsf m-1}, -1 \}$ is referred to as the regularity vector. A value of $-1$ for $\alpha_i$ indicates discontinuity of the basis functions at $\zeta_i$. We introduce the function space $\mathcal R^{\mathsf p}_{\bm \alpha} := \textup{span}\{\mathsf R_{i}^{\mathsf p}\}_{i=1}^n$, where the notation $\mathcal R^{\mathsf p}_{\alpha}$ is used to indicate that $\alpha_i = \alpha$ for $i=2,\cdots, \mathsf m-1$, suggesting continuity $C^{\alpha}$ for the spline function spaces. The construction of multivariate B-spline and NURBS basis functions follows a tensor-product manner. For $l=1, 2, \cdots, d$, given $\mathsf p_{l}$, $\mathsf n_{l}$, the knot vectors $\Xi_{l} = \{\xi_{1,l}, \xi_{2,l}, \cdots , \xi_{\mathsf n_{l}+\mathsf p_{l}+1, l}\}$, and the weight vectors $\{\mathsf w_{1,l}, \mathsf w_{2,l} \cdots , \mathsf w_{\mathsf n_{l},l}\}$, the univariate NURBS basis functions $\mathsf R^{\mathsf p_l}_{i_l,l}$ are well-defined. Consequently, the multivariate NURBS basis functions can be defined by exploiting the tensor product structure,
\begin{align*}
\mathsf R^{\mathsf p_1, \mathsf p_2\cdots , \mathsf p_d}_{i_1, i_2, \cdots ,i_d}(\xi_1,\cdots ,\xi_d) := \mathsf R^{\mathsf p_1}_{i_1, 1}(\xi_1) \otimes \mathsf R^{\mathsf p_2}_{i_2, 2}(\xi_2) \otimes \cdots \otimes \mathsf R^{\mathsf p_{d}}_{i_{d}, d}(\xi_d), \mbox{ for } i_l = 1,2, \cdots, \mathsf n_l \mbox{ and } l = 1, 2, \dots, d.
\end{align*}
The tensor product NURBS space is denoted as
\begin{align*}
\mathcal R^{\mathsf p_1, \mathsf p_2, \cdots \mathsf p_d}_{\bm \alpha_1, \bm \alpha_2, \cdots ,\bm \alpha_d} := \mathcal R^{\mathsf p_1}_{\bm \alpha_1} \otimes \mathcal R^{\mathsf p_2}_{\bm \alpha_2} \otimes \cdots \otimes \mathcal R^{\mathsf p_d}_{\bm \alpha_d} = \textup{span}\{ \mathsf R^{\mathsf p_1, \mathsf p_2, \cdots ,\mathsf p_d}_{i_1, i_2, \cdots , i_d} \}_{i_1=1, i_2=1, \cdots , i_d =1}^{\mathsf n_1, n_2, \cdots , \mathsf n_d}.
\end{align*}

\subsubsection{Mixed Formulation}
In this work, we always consider three-dimensional problems (i.e. $d=3$). Two discrete function spaces $\hat{\mathcal V}_h$ and $\hat{\mathcal Q}_h$ can be defined on the parametric domain $\hat{\Omega} = (0,1)^3$ as
\begin{align*}
\hat{\mathcal V}_h :=&
\mathcal R^{\mathsf p+1, \mathsf p+1, \mathsf p+1}_{\bm \alpha_1, \bm \alpha_2, \bm \alpha_3} \times \mathcal R^{\mathsf p+1, \mathsf p+1, \mathsf p+1}_{\bm \alpha_1, \bm \alpha_2, \bm \alpha_3} \times \mathcal R^{\mathsf p+1, \mathsf p+1, \mathsf p+1}_{\bm \alpha_1, \bm \alpha_2, \bm \alpha_3}, \quad
\hat{\mathcal Q}_h := \mathcal R^{\mathsf p, \mathsf p, \mathsf p}_{\bm \alpha_1, \bm \alpha_2, \bm \alpha_3}.
\end{align*}
Assuming the physical domain can be parametrized by a geometrical mapping $\bm \psi : \hat{\Omega} \rightarrow \Omega_{\bm X}$, the discrete functions on the physical domain can be defined through the pull-back operations,
\begin{align*}
\mathcal V_h :=  \{ \bm w : \bm w \circ \bm \psi \in \hat{\mathcal V}_h \}, \quad \mathcal Q_h := \{ q : q \circ \bm \psi \in \hat{\mathcal Q}_h \}.
\end{align*}
This pair of elements can be viewed as a smooth generalization of the Taylor-Hood element. We may then define the trial solution spaces for the velocity and pressure as
\begin{align*}
\mathcal S_{\bm v_h} &= \left\lbrace \bm v_h : \bm v_h(\cdot, t) \in \mathcal V_h, \bm v_h(\cdot,t) = \bm g \mbox{ on } \Gamma_{g}, t \in [0,T] \right\rbrace , \quad \mathcal S_{p_h} = \left\lbrace p_h : p_h(\cdot, t) \in \mathcal Q_h, t \in [0,T] \right\rbrace, 
\end{align*}
and the corresponding test function spaces are defined as
\begin{align*}
\mathcal V_{\bm v_h} &= \left\lbrace \bm w_{h} : \bm w_{h} \in \mathcal V_h, \bm w_{\bm v_h} = \bm 0 \mbox{ on } \Gamma_{g} \right\rbrace , \quad \mathcal V_{p_h} = \left\lbrace q_h : q_h \in \mathcal Q_h \right\rbrace.
\end{align*}
The semi-discrete formulation can be stated as follows. Find $\left\lbrace  \bm v_h(t), p_h(t) \right\rbrace^T \in \mathcal S_{\bm v_h} \times \mathcal S_{p_h}$ such that for $t\in [0, T]$,
\begin{align}
\label{eq:momentum_current}
& 0 = \mathbf B^m\left( \bm w_{h}; \dot{\bm v}_h, \bm v_h, p_h  \right) := \int_{\Omega} \bm w_{h} \cdot \rho \left( \frac{\partial \bm v_h}{\partial t} + \bm v_h \cdot \nabla \bm v_h - \bm f \right) + 2\mu \bm \varepsilon(\bm w_h) : \bm \varepsilon(\bm v_h)   - \nabla_{\bm x} \cdot \bm w_{h} \ p_h d\Omega_{\bm x} - \int_{\Gamma_{h}} \bm w_{h} \cdot \bm h d\Gamma_{\bm x}, \displaybreak[2] \\
\label{eq:mass_current}
& 0 = \mathbf B^c\left( q_h; \bm v_h \right) := \int_{\Omega} q_h \nabla_{\bm x} \cdot \bm v_h d\Omega_{\bm x},
\end{align}
for $\forall \left\lbrace \bm w_{h}, q_h \right\rbrace \in \mathcal V_{\bm v_h} \times \mathcal V_{p_h}$, with $\left\lbrace \bm v_{h}(0), p_{h}(0) \right\rbrace^T := \left\lbrace \bm v_{h0}, p_{h0} \right\rbrace^T$. Here $\bm v_{h0}$ and $p_{h0}$ are the $L_2$ projections of the initial data onto the finite dimensional trial solution spaces.

\subsection{The Generalized-$\alpha$ Schemes Under Comparison}
In this section, we present the two generalized-$\alpha$ schemes to be compared as time integration schemes for the semi-discrete formulation in \eqref{eq:momentum_current} and \eqref{eq:mass_current}. Let the time interval $[0,T]$ be divided into a set of $N_{\textup{ts}}$ subintervals delimited by a discrete time vector $\left\lbrace t_n \right\rbrace_{n=0}^{N_{\mathrm{ts}}}$. The time step size is defined as $\Delta t_n := t_{n+1} - t_n$. The discrete approximations to the velocity, its time derivative, and pressure at the time step $t_n$ are denoted as $\bm v_n$, $\dot{\bm v}_n$, and $p_n$. Let $N_A$ and $M_B$ represent basis functions for the velocity and pressure spaces, respectively, and let $\bm e_i$ be the Cartesian basis vector with $i=1,2,3$. We may define the residual vectors as follows,
\begin{align*}
\boldsymbol{\mathrm R}^m\left( \dot{\bm v}_{n}, \bm v_{n}, p_n \right) := \left\lbrace \mathbf B^m \left( N_A \bm e_i ;  \dot{\bm v}_{n}, \bm v_{n}, p_n \right) \right\rbrace, \quad \boldsymbol{\mathrm R}^c\left( \bm v_{n} \right) := \left\lbrace \mathbf B^c \left( M_B ;  \bm v_{n} \right) \right\rbrace.
\end{align*}
With these residual definitions, we introduce the two generalized-$\alpha$ schemes as follows. 

\paragraph{Scheme-1}
At time step $t_n$, given $\dot{\bm v}_n$, $\bm v_n$, $p_n$ and the time step size $\Delta t_n$, find $\dot{\bm v}_{n+1}$, $\bm v_{n+1}$, and $p_{n+1}$ such that,
\begin{align*}
& \boldsymbol{\mathrm R}^m\left( \dot{\bm v}_{n+\alpha_m}, \bm v_{n+\alpha_f}, p_{n+1} \right) = \bm 0, \quad \boldsymbol{\mathrm R}^c\left( \bm v_{n+\alpha_f} \right) = \bm 0, \quad \bm v_{n+1} = \bm v_{n} + \Delta t_n \dot{\bm v}_n + \gamma \Delta t_n \left( \dot{\bm v}_{n+1} - \dot{\bm v}_{n}\right), \\
& \dot{\bm v}_{n+\alpha_m} = \dot{\bm v}_{n} + \alpha_m \left(\dot{\bm v}_{n+1} - \dot{\bm v}_{n} \right), \quad \bm v_{n+\alpha_f} = \bm v_{n} + \alpha_f \left( \bm v_{n+1} - \bm v_{n} \right).
\end{align*}
Although the quantity $\dot{p}_{n+1}$ does not appear explicitly in the above formulation, we may still obtain a consistent definition as $\dot{p}_{n+1} := \left( p_{n+1} - p_n \right)/\Delta t_n$.

\paragraph{Scheme-2}
At time step $t_n$, given $\dot{\bm v}_n$, $\bm v_n$, $p_n$, and the time step size $\Delta t_n$, find $\dot{\bm v}_{n+1}$, $\bm v_{n+1}$, and $p_{n+1}$ such that,
\begin{align*}
& \boldsymbol{\mathrm R}^m\left( \dot{\bm v}_{n+\alpha_m}, \bm v_{n+\alpha_f}, p_{n+\alpha_f} \right) = \bm 0, \quad \boldsymbol{\mathrm R}^c\left( \bm v_{n+\alpha_f} \right) = \bm 0, \quad \bm v_{n+1} = \bm v_{n} + \Delta t_n \dot{\bm v}_n + \gamma \Delta t_n \left( \dot{\bm v}_{n+1} - \dot{\bm v}_{n}\right), \\
& \dot{\bm v}_{n+\alpha_m} = \dot{\bm v}_{n} + \alpha_m \left(\dot{\bm v}_{n+1} - \dot{\bm v}_{n} \right), \quad \bm v_{n+\alpha_f} = \bm v_{n} + \alpha_f \left( \bm v_{n+1} - \bm v_{n} \right), \quad p_{n+\alpha_f} = p_{n} + \alpha_f \left( p_{n+1} - p_{n} \right).
\end{align*}
In this scheme, given $\dot{p}_n$, $\dot{p}_{n+1}$ is given by
\begin{align*}
\dot{p}_{n+1} := \frac{p_{n+1} - p_n}{\gamma \Delta t_n} + \left(1 - \frac{1}{\gamma} \right) \dot{p}_n.
\end{align*}

\noindent The definitions of the above two schemes can be completed by assigning values to the parameters $\alpha_m$, $\alpha_f$, and $\gamma$. It has been demonstrated that, for linear problems, second-order accuracy, unconditionally stability, and optimal high frequency dissipation can be achieved via the parametrization
\begin{align*}
\alpha_m = \frac12 \left( \frac{3-\varrho_{\infty}}{1+\varrho_{\infty}} \right), \quad \alpha_f = \frac{1}{1+\varrho_{\infty}}, \quad \gamma = \frac{1}{2} + \alpha_m - \alpha_f,
\end{align*}
wherein $\varrho_{\infty} \in [0,1]$ is the spectral radius of the amplification matrix at the highest mode \cite{Jansen2000}. In this study, $\varrho_{\infty}$ 
is fixed to be $0.5$. In conventional CFD calculations, it is unnecessary to calculate $\dot{p}$. In hemodynamic simulations, however, both $\dot{\bm v}$ and $\dot{p}$ are needed as inputs to reduced-order models of the upstream and downstream vasculature, which are commonly used as boundary conditions coupled to the three-dimensional domain \cite{Moghadam2013}. We therefore examine the accuracy of $\dot{\bm v}$ and $\dot{p}$ in addition to $\bm v$ and $p$ in this work.

\section{Numerical Examples}
\label{sec:numerical_examples}
The parameters and results reported below are presented in the centimeter-gram-second units. 

\subsection{The Three-Dimensional Ethier-Steinman Benchmark}
\label{subsec:ES-benchmark}
In this section, we adopt the Ethier-Steinman solution \cite{Ethier1994} as a known exact solution to evaluate the accuracy of the two generalized-$\alpha$ schemes via the method of manufactured solutions. The analytic solutions take the following forms,
\begin{align*}
\bm v =& \begin{bmatrix}
-a \left( e^{ax}\sin(ay+dz) + e^{az} \cos(ax+dy) \right) \\
-a \left( e^{ay}\sin(az+dx) + e^{ax} \cos(ay+dz) \right) \\
-a \left( e^{az}\sin(ax+dy) + e^{ay} \cos(az+dx) \right)
\end{bmatrix} e^{-\nu d^2 t}, \displaybreak[2] \\
p =& -\frac{a^2}{2}\left( e^{2ax} + e^{2ay} + e^{2az} + 2\sin(ax+dy)\cos(az+dx)e^{a(y+z)} \right. \\
& \left. + 2\sin(ay+dz)\cos(ax+dy)e^{a(z+x)} + 2\sin(az+dx)\cos(ay+dz)e^{a(x+y)}  \right) e^{-2\nu d^2 t}.
\end{align*}
Consistent with the original design of the benchmark, the domain is a cube $\Omega = (-1,1)^3$. The constants are fixed as $a=\pi/4$ and $d=\pi/2$, and the fluid density $\rho$ and viscosity $\mu$ are fixed to be $1.0$ and $0.1$. Traction boundary conditions are enforced on all faces, and initial conditions are obtained from the exact solution. Prior to examination of the temporal accuracy, a spatial convergence study was first performed, and the code was verified to yield optimal convergence rates. In the following, we report the errors for a fixed mesh with $25$ elements along each dimension, polynomial degree $\mathsf p = 4$, and continuity $\alpha = 3$. The problem is simulated up to $T=1.0$ with various uniform time step sizes. The relative errors of the velocity, pressure, and their time derivatives in the $L_2$- and $H_1$-norms at time $t=1.0$ for Scheme-1 and Scheme-2 are presented in Tables \ref{table:temporal_ES_NS_3D_scheme_1} and \ref{table:temporal_ES_NS_3D_scheme_2}, respectively. Based on the numerical results, we make the following salient observations. First, the pressure converges \textit{linearly} with respect to the time step size in Scheme-1 and \textit{quadratically} in Scheme-2. Second, both schemes exhibit first-order temporal accuracy for the time derivative of pressure $\dot{p}_h$, yet Scheme-2 yields smaller errors for $\dot{p}_h$. Third, the velocity errors are identical in the two schemes, signifying that the degraded accuracy for pressure does not affect the accuracy for velocity. This phenomenon may, however, be attributed to the use of inf-sup stable elements. We do not rule out the possibility for velocity errors to be polluted by the use of Scheme-1 in \textit{stabilized} formulations.

\begin{table}[htbp]
\footnotesize
\begin{center}
\tabcolsep=0.21cm
\renewcommand{\arraystretch}{1.2}
\begin{tabular}{c | c c c c c c}
\hline
\hline
$N_{ts}$ & $10$ & $20$  & $40$ & $50$ & $80$ & $100$ \\
\hline
$\|\bm v - \bm v_h\|_{L_2(\Omega)} / \|\bm v\|_{L_2(\Omega)} $ & $2.85 \times 10^{-4}$ & $7.22 \times 10^{-5}$ & $1.82 \times 10^{-5}$ & $1.17\times 10^{-5}$ & $4.57 \times 10^{-6}$ & $2.92 \times 10^{-6}$ \\
order & - & 1.98 & 1.99 & 1.98 & 2.00 & 2.01 \\ 
$\|\bm v - \bm v_h\|_{H_1(\Omega)} / \|\bm v\|_{H_1(\Omega)} $ & $2.44 \times 10^{-4}$ & $6.19 \times 10^{-5}$ & $1.56 \times 10^{-5}$ & $9.99\times 10^{-6}$ & $3.91 \times 10^{-6}$ & $2.51 \times 10^{-6}$ \\
order & - & 1.98 & 1.99 & 2.00 & 2.00 & 1.99 \\ 
$\| p - p_h\|_{L_2(\Omega)} / \| p \|_{L_2(\Omega)} $ & $1.66 \times 10^{-2}$ & $8.27 \times 10^{-3}$ & $4.12 \times 10^{-3}$ & $3.30 \times 10^{-3}$ & $2.06 \times 10^{-3}$  & $1.65 \times 10^{-3}$  \\
order & - & 1.01 & 1.01 & 0.99 & 1.00 & 0.99 \\
$\| p - p_h\|_{H_1(\Omega)} / \| p \|_{H_1(\Omega)} $ & $1.66 \times 10^{-2}$ & $8.26 \times 10^{-3}$ & $4.12 \times 10^{-3}$ & $3.30 \times 10^{-3}$ & $2.06 \times 10^{-3}$ & $1.65 \times 10^{-3}$ \\
order & - & 1.01 & 1.00 & 0.99 & 1.00 & 0.99 \\ 
\hline
$\|\dot{\bm v} - \dot{\bm v}_h\|_{L_2(\Omega)} / \|\dot{\bm v}\|_{L_2(\Omega)} $ &  $3.66 \times 10^{-3}$ & $1.92 \times 10^{-3}$ & $9.90 \times 10^{-4}$ & $7.98 \times 10^{-4}$ & $5.04 \times 10^{-4}$ & $4.05 \times 10^{-4}$ \\
order & - & 0.93 & 0.96 & 0.97 & 0.98 & 0.98 \\
$\|\dot{\bm v} - \dot{\bm v}_h\|_{H_1(\Omega)} / \|\dot{\bm v}\|_{H_1(\Omega)} $ &$4.20 \times 10^{-3}$ & $2.05 \times 10^{-3}$ & $1.02 \times 10^{-3}$ & $8.20 \times 10^{-4}$ & $5.13 \times 10^{-4}$ & $4.10 \times 10^{-4}$ \\
order & - & 1.03 & 1.01  & 0.98 & 1.00 & 1.00 \\
$\|\dot{p} - \dot{p}_h\|_{L_2(\Omega)} / \|\dot{p}\|_{L_2(\Omega)} $ & $4.20 \times 10^{-2}$ & $2.08 \times 10^{-2}$ & $1.03 \times 10^{-2}$ & $8.26 \times 10^{-3}$ & $5.15 \times 10^{-3}$  & $4.12 \times 10^{-3}$  \\
order & - & 1.01 & 1.01 & 0.99 & 1.01 & 1.00 \\
$\|\dot{p} - \dot{p}_h\|_{H_1(\Omega)} / \|\dot{p}\|_{H_1(\Omega)} $ & $4.20 \times 10^{-2}$ & $2.08 \times 10^{-2}$ & $1.03 \times 10^{-2}$ & $8.26 \times 10^{-3}$ & $5.15 \times 10^{-3}$  & $4.12 \times 10^{-3}$  \\
order & - & 1.01 & 1.01 & 0.99 & 1.01 & 1.00 \\ 
\hline 
\hline
\end{tabular}
\end{center}
\caption{Temporal convergence rates for the Ethier-Steinman benchmark using Scheme-1.}
\label{table:temporal_ES_NS_3D_scheme_1}
\end{table}

\begin{table}[htbp]
\footnotesize
\begin{center}
\tabcolsep=0.21cm
\renewcommand{\arraystretch}{1.2}
\begin{tabular}{c | c c c c c c}
\hline
\hline
$N_{ts}$ & $10$ & $20$  & $40$ & $50$ & $80$ & $100$ \\
\hline
$\|\bm v - \bm v_h\|_{L_2(\Omega)} / \|\bm v\|_{L_2(\Omega)} $ & $2.85 \times 10^{-4}$ & $7.22 \times 10^{-5}$ & $1.82 \times 10^{-5}$ & $1.17\times 10^{-5}$ & $4.57 \times 10^{-6}$ & $2.92 \times 10^{-6}$ \\
order & - & 1.98 & 1.99 & 1.98 & 2.00 & 2.01 \\ 
$\|\bm v - \bm v_h\|_{H_1(\Omega)} / \|\bm v\|_{H_1(\Omega)} $ & $2.44 \times 10^{-4}$ & $6.19 \times 10^{-5}$ & $1.56 \times 10^{-5}$ & $9.99\times 10^{-6}$ & $3.91 \times 10^{-6}$ & $2.51 \times 10^{-6}$ \\
order & - & 1.98 & 1.99 & 2.00 & 2.00 & 1.99 \\ 
$\| p - p_h\|_{L_2(\Omega)} / \| p \|_{L_2(\Omega)} $ & $2.39 \times 10^{-4}$ & $5.98 \times 10^{-5}$ & $1.49 \times 10^{-5}$ & $9.54\times 10^{-6}$ & $3.72 \times 10^{-6}$ & $2.38 \times 10^{-6}$ \\
order & - & 2.00 & 2.00 & 2.00 & 2.00 & 2.00 \\
$\| p - p_h\|_{H_1(\Omega)} / \| p \|_{H_1(\Omega)} $ & $2.63 \times 10^{-4}$ & $6.58 \times 10^{-5}$ & $1.64 \times 10^{-5}$ & $1.05\times 10^{-5}$ & $4.13 \times 10^{-6}$ & $2.65 \times 10^{-6}$ \\
order & - & 2.00 & 2.00 & 2.00 & 1.99 & 1.99 \\ 
\hline
$\|\dot{\bm v} - \dot{\bm v}_h\|_{L_2(\Omega)} / \|\dot{\bm v}\|_{L_2(\Omega)} $ &  $3.66 \times 10^{-3}$ & $1.92 \times 10^{-3}$ & $9.90 \times 10^{-4}$ & $7.98 \times 10^{-4}$ & $5.04 \times 10^{-4}$ & $4.05 \times 10^{-4}$ \\
order & - & 0.93 & 0.96 & 0.97 & 0.98 & 0.98 \\
$\|\dot{\bm v} - \dot{\bm v}_h\|_{H_1(\Omega)} / \|\dot{\bm v}\|_{H_1(\Omega)} $ &$4.20 \times 10^{-3}$ & $2.05 \times 10^{-3}$ & $1.02 \times 10^{-3}$ & $8.20 \times 10^{-4}$ & $5.13 \times 10^{-4}$ & $4.10 \times 10^{-4}$ \\
order & - & 1.03 & 1.01  & 0.98 & 1.00 & 1.00 \\
$\|\dot{p} - \dot{p}_h\|_{L_2(\Omega)} / \|\dot{p}\|_{L_2(\Omega)} $ & $7.33 \times 10^{-3}$ & $3.99 \times 10^{-3}$ & $2.03 \times 10^{-3}$ & $1.63 \times 10^{-3}$ & $1.02 \times 10^{-3}$  & $8.18 \times 10^{-4}$ \\
order & - & 0.88 & 0.97 & 0.98 & 1.00 & 0.99 \\
$\|\dot{p} - \dot{p}_h\|_{H_1(\Omega)} / \|\dot{p}\|_{H_1(\Omega)} $ & $7.41 \times 10^{-3}$ & $4.00 \times 10^{-3}$ & $2.03 \times 10^{-3}$ & $1.63 \times 10^{-3}$ & $1.02 \times 10^{-3}$  & $8.18 \times 10^{-4}$  \\
order & - & 0.89 & 0.98 & 0.98 & 1.00 & 0.99 \\ 
\hline 
\hline
\end{tabular}
\end{center}
\caption{Temporal convergence rates for the Ethier-Steinman benchmark using Scheme-2.}
\label{table:temporal_ES_NS_3D_scheme_2}
\end{table}

\begin{figure}
	\begin{center}
	\begin{tabular}{c}
\includegraphics[angle=0, trim=100 80 100 120, clip=true, scale = 0.3]{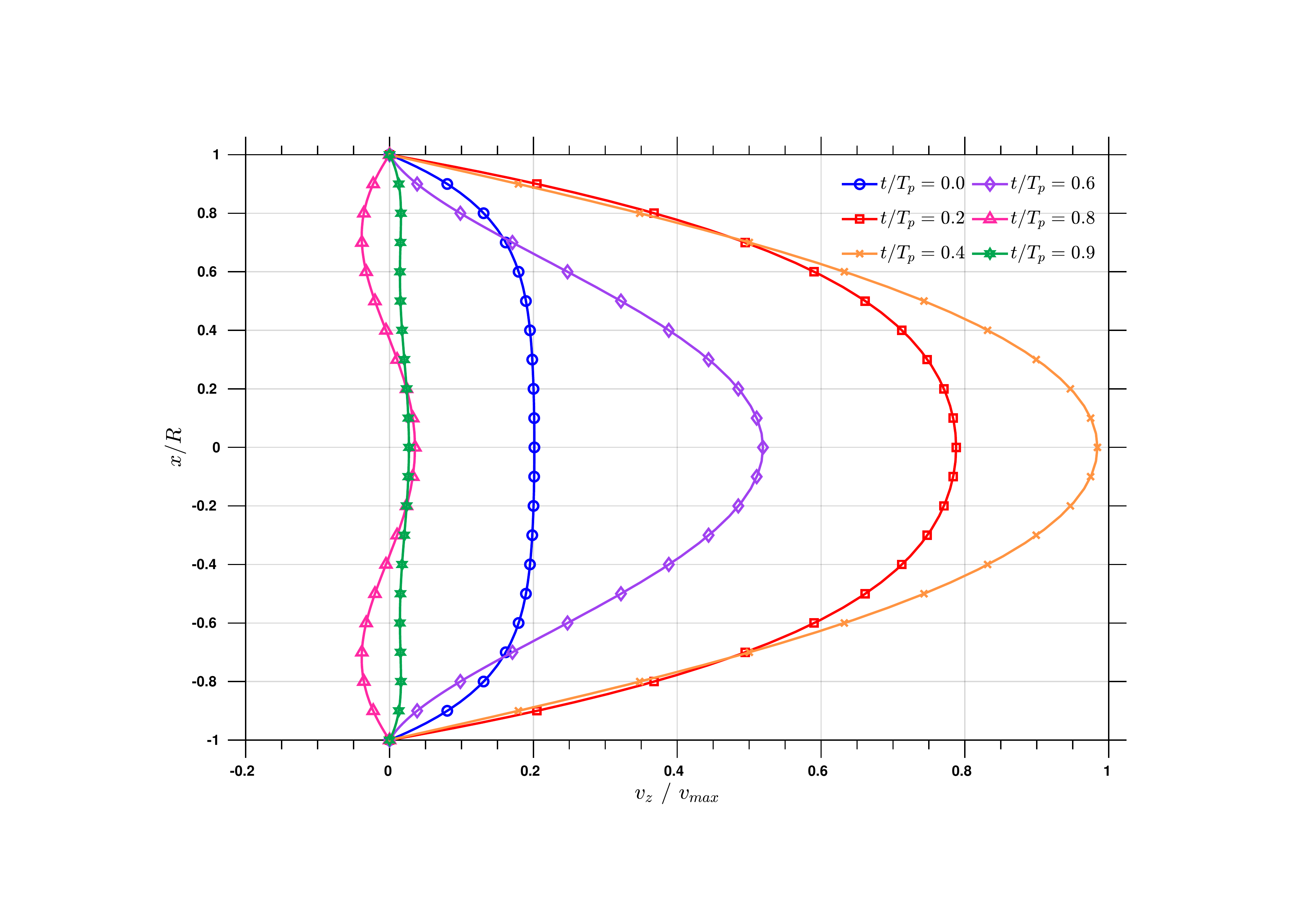} 
\end{tabular}
\caption{Axial velocity profiles along the $x$-axis at $z=0.5$ of the pipe at different time instances. All axial velocities are normalized to the maximum axial velocity $v_{max}$.} 
\label{fig:womersley_velo}
\end{center}
\end{figure}

\subsection{Pulsatile Flow in a Rigid Pipe}
\label{subsec:pulastile_flow}
\begin{figure}
	\begin{center}
	\begin{tabular}{cc}
\includegraphics[angle=0, trim=90 85 90 80, clip=true, scale = 0.36]{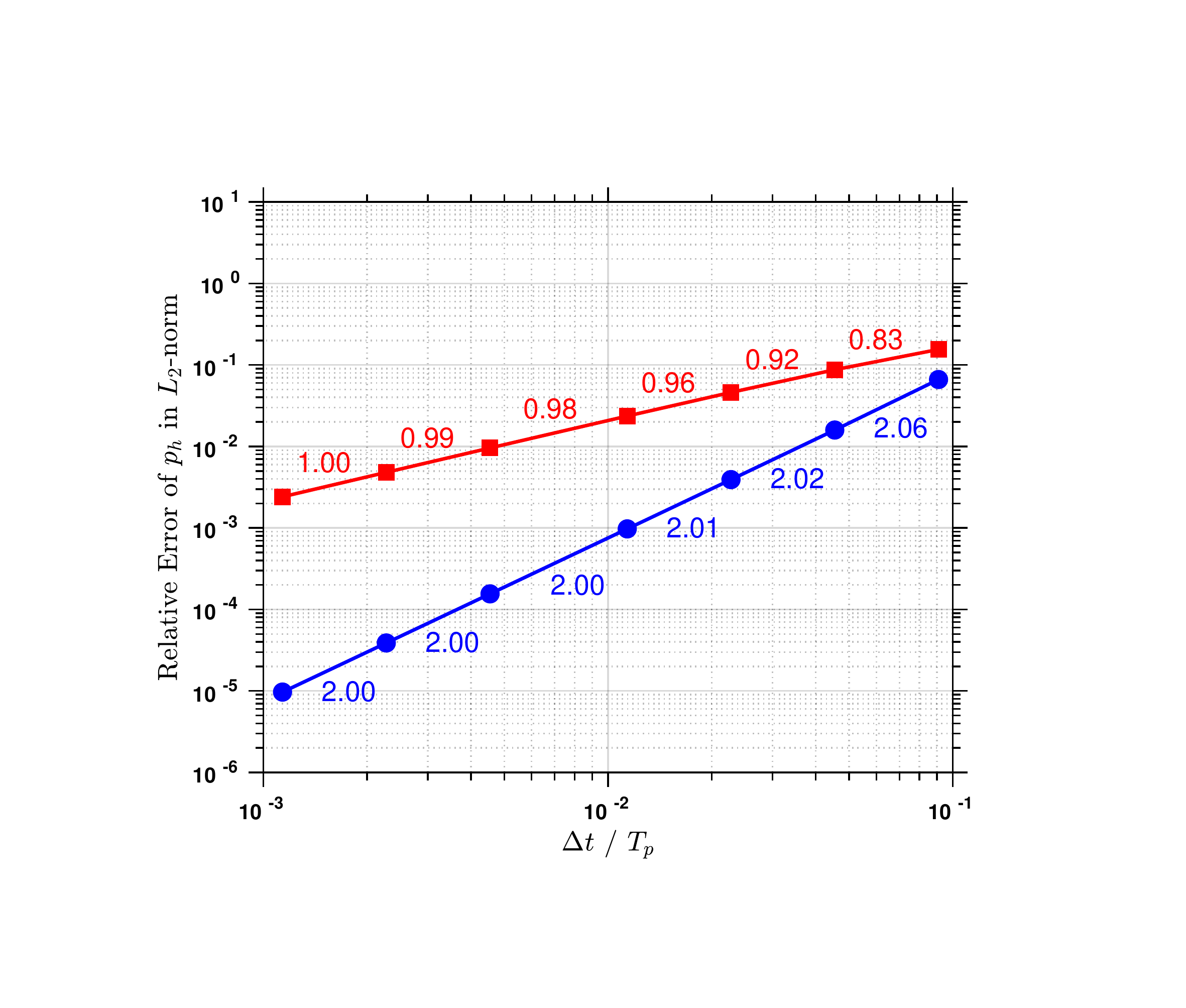} &
\includegraphics[angle=0, trim=90 85 90 80, clip=true, scale = 0.36]{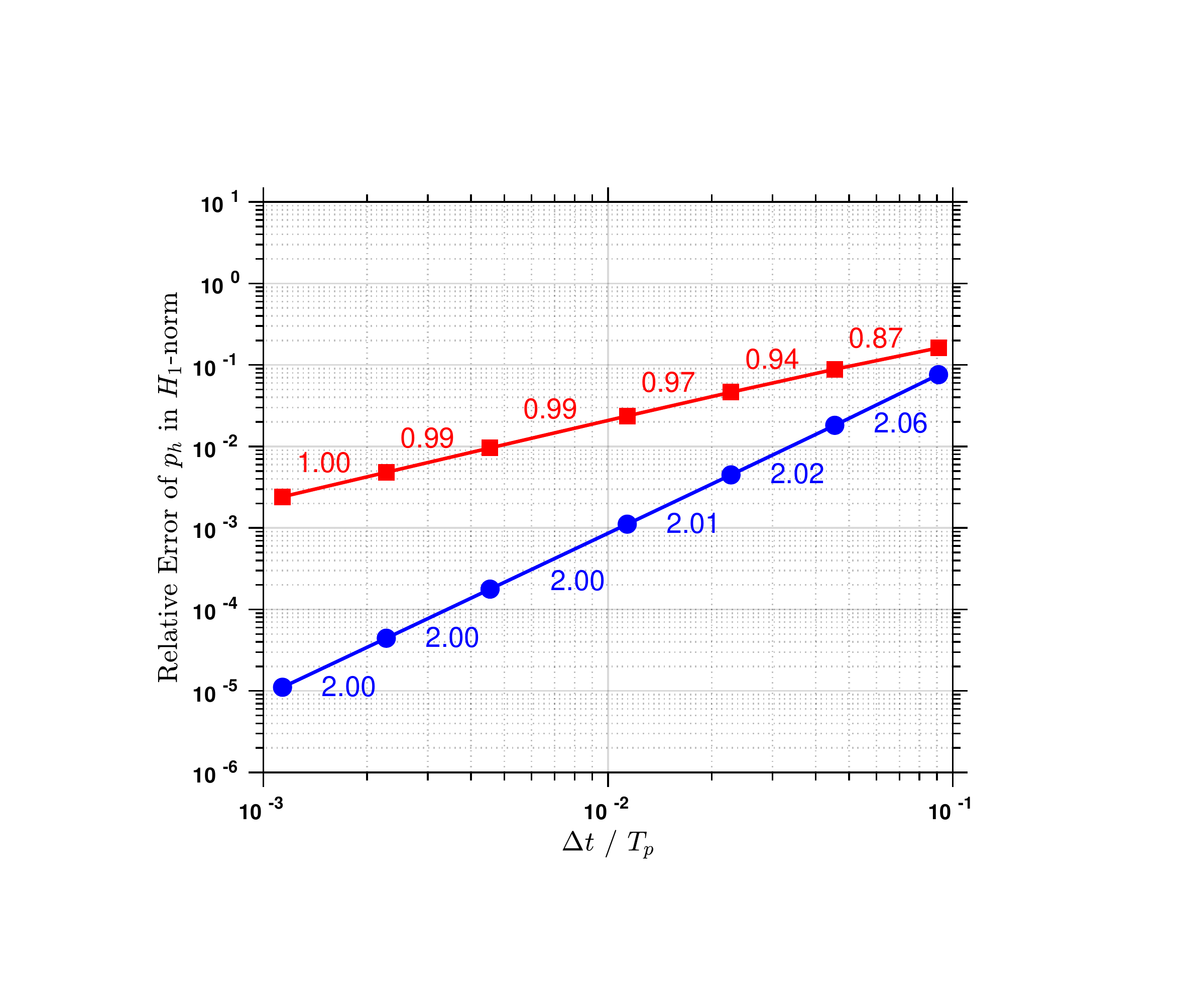} \\
(a)  & (b) \\
\includegraphics[angle=0, trim=90 85 90 80, clip=true, scale = 0.36]{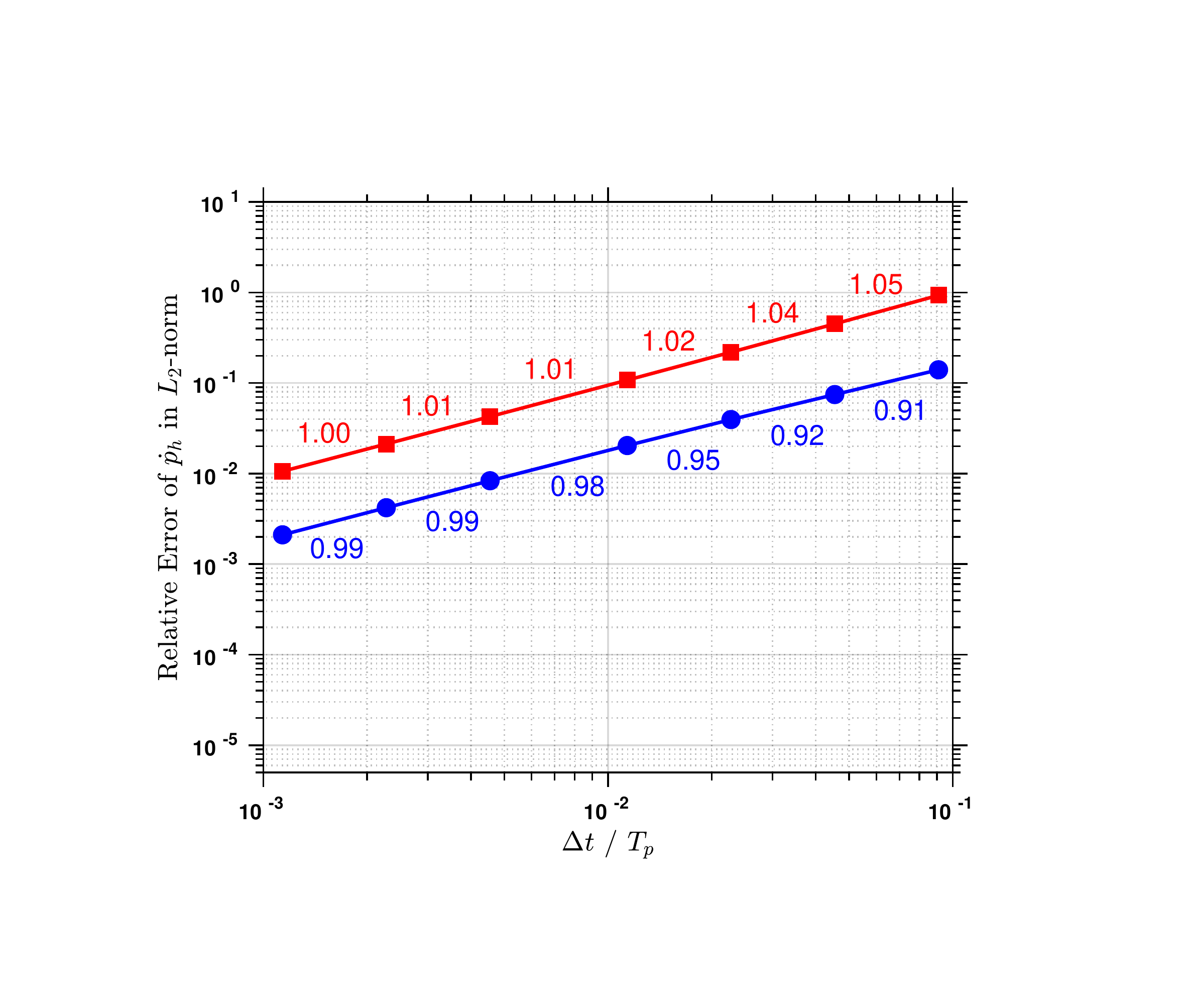} &
\includegraphics[angle=0, trim=90 85 90 80, clip=true, scale = 0.36]{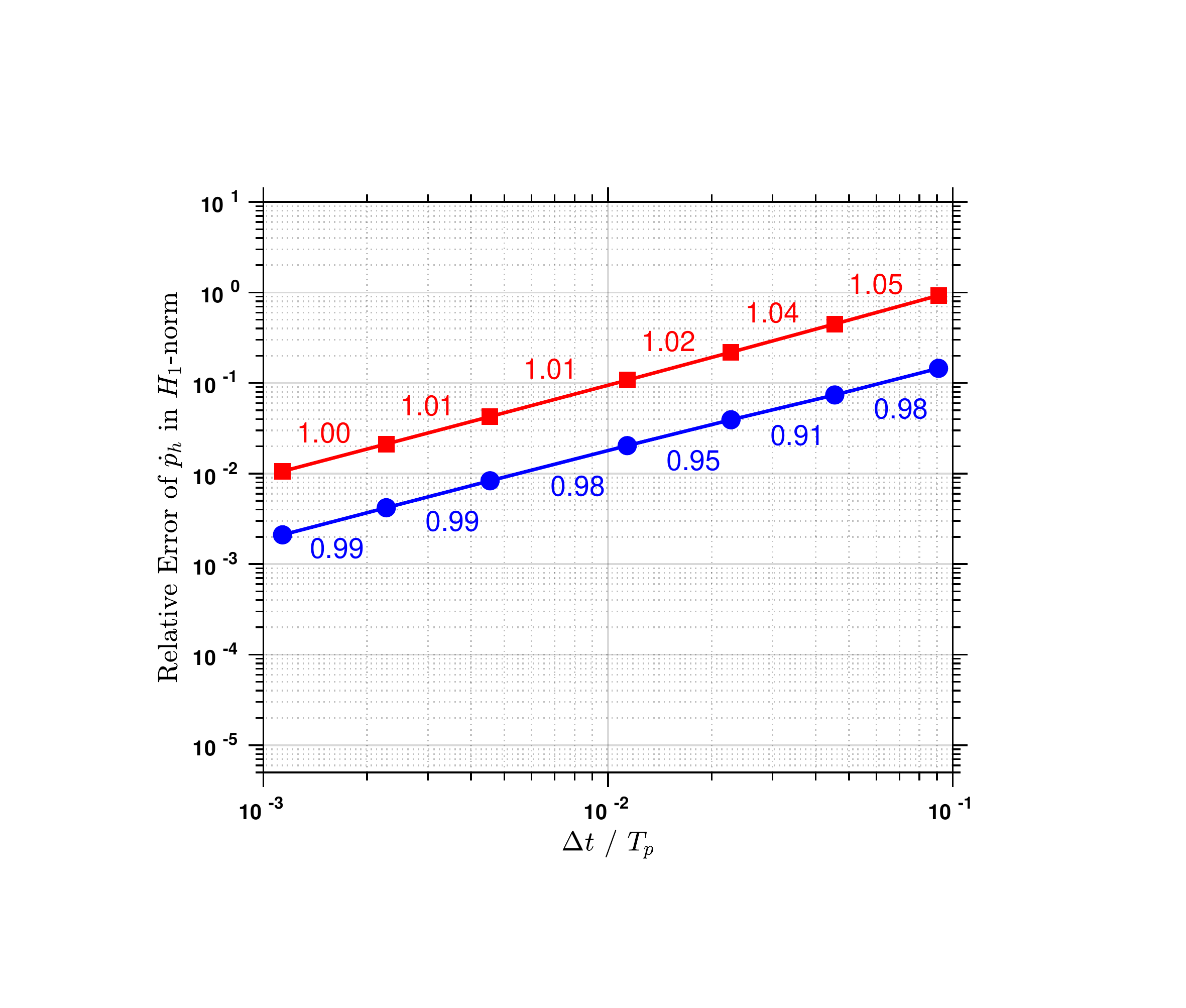} \\
(c)  & (d) \\
\multicolumn{2}{c}{ \includegraphics[angle=0, trim=550 110 550 150, clip=true, scale = 0.5]{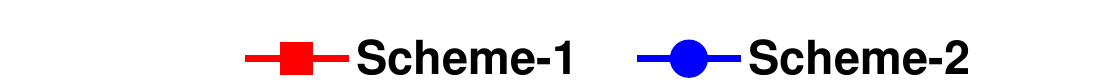} }
\end{tabular}
\caption{The relative errors of (a) $p_h$ in $L_2$ norm, (b) $p_h$ in $H_1$ norm, (c) $\dot{p}_h$ in $L_2$ norm, and (d) $\dot{p}_h$ in $H_1$ norm with different time step sizes. The convergence rates calculated from adjacent errors are annotated.} 
\label{fig:conv_pressure}
\end{center}
\end{figure}

In the second example, we utilize the Womersley solution to evaluate the two schemes for flows under more physiological settings. The Womersley solution describes axisymmetric, fully developed flow subject to a periodically oscillating pressure gradient in a rigid cylindrical pipe with circular cross section \cite{Womersley1955}. Here, we consider the $z$-axis to be oriented along the length of the pipe and express the pressure in the following Fourier series,
\begin{align}
\label{eq:womersley_pressure}
p = p_{\textup{ref}} + \left( k_0 + \sum_{n=1}^{N} k_n e^{\iota n\omega t} \right) z,
\end{align}
where $p_{\textup{ref}}$ is the reference pressure at $z=0$, $k_0$ is the steady component of the pressure gradient, $k_n$ is the amplitude of the oscillatory component of the pressure gradient at the $n$-th Fourier mode, $\iota$ is the solution of the equation $\iota^2=-1$, $\omega:=2\pi / T_p$ is the fundamental frequency of the pressure wave, and $T_p$ is the period of the pressure. Per the assumptions of axisymmetric and fully developed flow, the velocity is identically zero in the radial and circumferential directions and takes the following analytical form in the axial direction,
\begin{align}
\label{eq:womersley_velo}
v_z = \frac{k_0}{4\mu}\left(r^2 - R^2\right) + \sum\limits_{n=1}^{N} \frac{\iota k_n}{\rho n\omega} \left( 1 - \frac{J_0(\iota^{\frac{3}{2}} \alpha_n \frac{r}{R})}{J_0(\iota^{\frac{3}{2}}\alpha_n)} \right) e^{\iota n \omega t},
\end{align}
wherein $r:= \sqrt{x^2+y^2}$, $R$ is the radius of the pipe, $J_0$ is the zeroth-order Bessel function of the first kind, and $\alpha_n := R \sqrt{n\omega / \nu}$ is the Womersley number for the $n$-th Fourier mode. The complex form of the solutions indicates the existence of two sets of real independent solutions. We use the real parts of \eqref{eq:womersley_pressure} and \eqref{eq:womersley_velo} to serve as the benchmark solution. In this study, the pipe radius $R$ is $0.3$, and the pipe length is $1.0$. The fluid density $\rho$ and viscosity $\mu$ are fixed to be $1.0$ and $0.04$, respectively, and the period $T_p$ is set as $1.1$. The fundamental frequency $\omega$ and  the corresponding Womersley number $\alpha_1$ are approximately $5.71$ and $3.59$, respectively. The reference pressure $p_{\textup{ref}}$ is $0$. We choose to use a single mode in the Fourier series to represent the oscillatory part of the pressure (i.e. $N=1$), and the values of $k_0$ and $k_1$ are $-21.0469$ and $-33.0102+42.9332\iota$. The parameters above are chosen to represent physiological flows in arteries. The circular cross section is exactly represented by NURBS, and the parametrization is based on quadratic basis functions, which can be found in Example 6.2 of \cite{Takacs2011}. Using $k$-refinement, the spatial mesh for analysis consists of $10 \times 10$ elements in the $x$-$y$ plane and $16$ elements along the $z$ axis, with polynomial degree $\mathsf p = 4$ and continuity $\alpha = 3$. We apply the no-slip boundary condition on the wall and traction boundary conditions on both ends of the pipe. The analytical forms of the tractions as well as the initial conditions are obtained from the exact solution in \eqref{eq:womersley_pressure} and \eqref{eq:womersley_velo}. Simulations are performed with uniform time step sizes. Similar to the previous example, we observed that the velocity errors in all norms are identical for the two schemes. The velocity profiles at various time instances are illustrated in Figure \ref{fig:womersley_velo}. The relative errors of $p_h$ and $\dot{p}_{h}$ in the $L_2$- and $H_1$-norms at time $t=0.8$  are plotted in Figure \ref{fig:conv_pressure}. We observe that the pressure converges \textit{linearly} in Scheme-1 and \textit{quadratically} in Scheme-2 in both norms. While both schemes exhibit first-order temporal accuracy for $\dot{p}_{h}$, we again find that $\dot{p}_h$ calculated from Scheme-2 is more accurate.


\section{Concluding Remarks}
In this work, we implemented and compared two generalized-$\alpha$ schemes for the incompressible Navier-Stokes equations using inf-sup stable higher-order NURBS elements for spatial discretization. In the widely adopted scheme that is used in several existing CFD and FSI implementations, we observed only first-order accuracy for pressure. Although we have not observed degraded accuracy for velocity in these schemes, we cannot rule out this possibility for other spatial discretization methods, such as the stabilized finite element method. By evaluating the pressure at the intermediate time step in the temporal scheme, second-order temporal accuracy is recovered. Moreover, the time derivative of the pressure is computed more accurately. Therefore, we recommend the use of Scheme-2 as the generalized-$\alpha$ scheme of choice for integrating the incompressible Navier-Stokes equations.

\section*{Acknowledgements}
This work is supported by the NIH under the award numbers 1R01HL121754, 1R01HL123689, R01EB01830204, the computational resources from the Stanford Research Computing Center, and the Extreme Science and Engineering Discovery Environment supported by the NSF grant ACI-1053575. Ingrid Lan is supported by the NSF Graduate Research Fellowship and the Stanford Graduate Fellowship in Science and Engineering.

\bibliographystyle{plain}  
\bibliography{note_ga}

\end{document}